\documentclass[12pt]{article}
\usepackage{amsmath,amssymb,theorem,amstext,amsgen,amsbsy,amsopn,amsfonts,graphicx,cases}
\usepackage{graphicx}
\usepackage{subfigure}
\usepackage{psfrag}
\usepackage{color}
\usepackage{enumerate}
\usepackage{epstopdf}
\usepackage{array}

\usepackage{latexsym}
\newtheorem{thm}{Theorem}[section]

\newtheorem{cor}[thm]{Corollary}
\newtheorem{lem}[thm]{Lemma}
\newtheorem{prop}[thm]{Proposition}
\theorembodyfont{\rmfamily}

\def\pf{\bigskip\noindent {\bf Proof.}}

\def\dfn#1{{\sl #1}}

\def\less{\setminus}

\def\pf{\bigskip\noindent {\emph{Proof.}}}

\def\qed{ \hfill\vrule height3pt width6pt depth2pt}

\def\pf{\bigskip\noindent {\bf Proof.  }}

\title{Planar anti-Ramsey numbers of  paths and cycles}

\author{Yongxin Lan$^1$, Yongtang Shi$^1$, Zi-Xia Song$^2$\thanks{Corresponding author.}\\
$^1$Center for Combinatorics and LPMC\\
Nankai University, Tianjin 300071, China\\
 Email: lan@mail.nankai.edu.cn; shi@nankai.edu.cn\\
$^2$Department  of Mathematics\\
 University of Central Florida, Orlando, FL 32816, USA\\
   Email: Zixia.Song@ucf.edu
}
\date{August 31, 2017}
\begin{document}
\maketitle
\begin{abstract}
Motivated by anti-Ramsey numbers introduced by Erd\H{o}s, Simonovits and S\'os in 1975,  we study the anti-Ramsey problem when host graphs are plane triangulations. Given a positive integer $n$ and  a planar graph  $H$, let  $\mathcal{T}_n(H)$  be the family of  all plane triangulations $T$ on  $n$ vertices such that $T$  contains  a subgraph isomorphic to  $H$.  The \dfn{planar anti-Ramsey number of $H$}, denoted  $ar_{_\mathcal{P}}(n, H)$,  is the  maximum number  of colors in an edge-coloring of  a plane triangulation $T\in \mathcal{T}_n(H)$ such that $T$ contains   no    rainbow copy of $H$. Analogous to  anti-Ramsey numbers and Tur\'an numbers, planar anti-Ramsey numbers are closely related to  planar Tur\'an numbers, where the \dfn{planar Tur\'an number of $H$} is  the maximum number of edges of  a   planar graph on $n$ vertices without containing $H$ as a subgraph.  The  study of $ar_{_\mathcal{P}}(n, H)$  (under the name of rainbow numbers) was initiated by Hor\v{n}\'ak,  Jendrol$'$,  Schiermeyer and  Sot\'ak [J Graph Theory 78 (2015) 248--257]. In this paper we study planar anti-Ramsey numbers for paths and cycles. We first establish  lower bounds for  $ar_{_\mathcal{P}}(n, P_k)$ when  $n\ge k\ge8$.  We then improve the  existing lower bound  for $ar_{_\mathcal{P}}(n, C_k)$ when $k\geq 5$ and $n\geq k^2-k$.   Finally, using the main ideas in the above-mentioned paper, we obtain   upper bounds for $ar_{_\mathcal{P}}(n, C_6)$  when $n\ge8$ and $ar_{_\mathcal{P}}(n, C_7)$ when $n\geq 13$, respectively.  \end{abstract}


\baselineskip 17pt
\section{Introduction}

All graphs considered in this paper are finite and simple. Motivated by anti-Ramsey numbers introduced by Erd\H{o}s, Simonovits  and S\'os~\cite{ESS} in 1975,  we study the anti-Ramsey problem when host graphs are plane triangulations.  A subgraph of an edge-colored graph is \dfn{rainbow} if all of its edges have different colors.  Let $\mathcal{F}$ be a family of planar graphs. For the purpose of this paper, we call an edge-coloring that contains no rainbow copy  of any graph in $\mathcal{F}$ an \dfn{$\mathcal{F}$-free edge-coloring}. A  graph $G$ is \dfn{$\mathcal{F}$-free} if no subgraph of $G$ is isomorphic to any  graph in $\mathcal  {F}$.  Let  $n_{_\mathcal{F}}$ be the smallest integer such that for any $n\ge  n_{_\mathcal{F}}$, there exists a  plane triangulation on  $n$ vertices that is not $\mathcal{F}$-free. Such an integer $n_{_\mathcal{F}}$ is well-defined, because for any $F\in \mathcal{F}$, we can obtain a plane  triangulation from a plane drawing of $F$ by adding new edges.  When $\mathcal{F}=\{F\}$, then $n_{_\mathcal{F}} =|F|$. 
For any  integer $n\ge n_{_\mathcal{F}}$, let  $\mathcal{T}_n(\mathcal{F})$  be the family of  all plane triangulations $T$ on  $n$ vertices such that $T$  is not $\mathcal{F}$-free.  The \dfn{planar anti-Ramsey number of $\mathcal{F}$}, denoted
 $ar_{_\mathcal{P}}(n,\mathcal{F})$,  is the  maximum number  of colors in an $\mathcal{F}$-free edge-coloring of a plane triangulation in $\mathcal{T}_n(\mathcal{F})$.  Clearly,  $ar_{_\mathcal{P}}(n,\mathcal{F}) < 3n-6$.  It is worth noting that this problem becomes  trivial if the host plane triangulation  on $n$ vertices is  $\mathcal{F}$-free, because $3n-6$ colors can be used.  \medskip
 
 Analogous to the relation between anti-Ramsey numbers and Tur\'an numbers proved in~\cite{ESS}, planar anti-Ramsey numbers are closely related to planar Tur\'an numbers~\cite{Dowden}.  The \dfn{planar Tur\'an number of $\mathcal{F}$}, denoted  $ex_{_\mathcal{P}}(n,\mathcal{F})$, is  the maximum number of edges of  an  $\mathcal{F}$-free  planar graph on $n$ vertices. Given an edge-coloring $c$ of a host graph $T$ in $\mathcal{T}_n(\mathcal{F})$, we define a \dfn{representing graph} of $c$ to be a spanning subgraph $R$ of $T$ obtained by taking one edge of each color  under the coloring $c$ (where $R$ may contain isolated vertices). Clearly, if $c$ is an $\mathcal{F}$-free edge-coloring of $T$, then $R$ is $\mathcal{F}$-free. Thus $ar_{_\mathcal{P}}(n,\mathcal{F})\le ex_{_\mathcal{P}}(n,\mathcal{F})$ for any $n\ge n_{_\mathcal{F}}$. When $\mathcal{F}$ consists of a single graph $H$, we write  $ar_{_\mathcal{P}}(n, H)$ and $ex_{_\mathcal{P}}(n, H)$ instead of  $ar_{_\mathcal{P}}(n,\{H\})$ and $ex_{_\mathcal{P}}(n, \{H\})$. Given a planar graph $H$, let $\mathcal{H}=\{H-e: \, e\in E(H)\}$.  Let $G$ be an $\mathcal{H}$-free plane subgraph of a plane triangulation $T\in \mathcal{T}_n(H)$  with $e(G)=ex_{_\mathcal{P}}(n, \mathcal{H})$.  We then obtain an $H$-free edge-coloring of $T$ by coloring the edges of $G$ with distinct colors and then coloring  the edges in $E(T)\less E(G)$ with a new color. Hence,  $1+ex_{_\mathcal{P}}(n,\mathcal{H})\le ar_{_\mathcal{P}}(n,H) $ for any $n\ge |H|$.  We obtain the following analogous result.\medskip
 
 \begin{prop}\label{LU} Given a planar graph $H$ and a positive integer $n\ge |H|$, 
 $$1+ex_{_\mathcal{P}}(n,\mathcal{H})\le ar_{_\mathcal{P}}(n,H) \le ex_{_\mathcal{P}}(n, H),$$
 where $\mathcal{H}=\{H-e: \, e\in E(H)\}$. 
 \end{prop}
 
 Colorings of plane graphs  that avoid rainbow faces have also been  studied,  see, e.g., \cite{Kral, JMSS, West, Zykov}. Various results on anti-Ramsey numbers can be found in:  \cite{Alon, Jiang2004, JST, Jiang2002, Jiang2009, Li, M, Ingo} to name a few. 
The  study of planar anti-Ramsey numbers $ar_{_\mathcal{P}}(n, H)$  was initiated by Hor\v{n}\'ak,  Jendrol$'$,  Schiermeyer and  Sot\'ak~\cite{HJSS} (under the name of rainbow numbers). We summarize  their results in~\cite{HJSS} as follows, where given two positive integers $a$ and $b$, we use  \dfn{$a\text{ mod }b$} to denote the remainder when  $a$ is divided by $b$.  We use $P_k$ and $C_k$ to denote the path and cycle on $k$ vertices, respectively.

\begin{thm}[\cite{HJSS}]\label{HJSS2015}  Let $n, k$ be positive integers. 
\begin{enumerate}[(a)]
\item  $ar_{_\mathcal{P}}(n, C_3)=\lfloor{(3n-6)/}2\rfloor$ for $n\geq 4$.\vspace{-8pt}
\item $ar_{_\mathcal{P}}(n, C_4)\leq 2(n-2)$ for $n\geq 4$, and $ar_{_\mathcal{P}}(n, C_4)\geq (9(n-2)-4r)/5 $ for $n\geq 42$ and $r=(n-2)\text{ mod }20$.\vspace{-8pt}
\item  $ar_{_\mathcal{P}}(n, C_5)\leq {5(n-2)}/2$ for $n\geq 5$, and $ar_{_\mathcal{P}}(n, C_5)\geq (19(n-2)-10r)/9$ for $n\geq 20$ and $r=(n-2)\text{ mod }18$.\vspace{-8pt}
\item  $ar_{_\mathcal{P}}(n, C_k)\geq (3n-6)\cdot \frac{k-3}{k-2}-\frac{2k-7}{k-2}$ for $6\leq k\leq n$.
\end{enumerate}
\end{thm}

Finding  exact values of $ar_{_\mathcal{P}}(n, H)$ is far from trivial. As observed in~\cite{HJSS}, an induction argument in general  cannot be applied to compute $ar_{_\mathcal{P}}(n, H)$ because  deleting a vertex from a plane triangulation may result in a graph that is no longer a plane triangulation. \medskip

 Dowden~\cite{Dowden} began the study of planar Tur\'an numbers $ex_{_\mathcal{P}}(n, H)$ (under the name of ``extremal" planar graphs) and  proved Theorem~\ref{Dowden2015} below, where each bound is tight. 
 \begin{thm}[\cite{Dowden}]\label{Dowden2015}  Let $n$ be a positive integer. 
\begin{enumerate}[(a)]
\item  $ex_{_\mathcal{P}}(n, C_3)=2n-4$ for $n\geq 3$. \vspace{-8pt}
\item $ex_{_\mathcal{P}}(n, C_4)\leq {15}(n-2)/7$ for $n\geq 4$. \vspace{-8pt}
\item  $ex_{_\mathcal{P}}(n, C_5)\leq (12n-33)/{5}$ for $n\geq 11$.
\end{enumerate}
\end{thm}

 By Proposition~\ref{LU} and Theorem~\ref{Dowden2015}(c), we see that  $ar_{_\mathcal{P}}(n, C_5)\le (12n-33)/{5}$ for $n\geq 11$. This improves the upper bound for $ar_{_\mathcal{P}}(n, C_5)$ in Theorem~\ref{HJSS2015}(c) when $n\ge11$.  Notice that the upper bound in Proposition~\ref{LU} in general  is quite loose,  for example,   $ex_{_\mathcal{P}}(n, C_3) -ar_{_\mathcal{P}}(n, C_3)=\lceil n/2\rceil-1$ for all $n\ge4$. In this paper we study planar anti-Ramsey numbers for paths and cycles. 
In Section~\ref{path}, we establish    lower bounds for  $ar_{_\mathcal{P}}(n, P_k)$ when  $n\ge k\ge8$.  
In  Section~\ref{cycle}, we first  improve the  existing lower bounds  for $ar_{_\mathcal{P}}(n, C_k)$ when $k\geq 5$ and $n\geq k^2-k$, which improves  Theorem~\ref{HJSS2015}(c,d).   We then use the main ideas in~\cite{HJSS} by studying lower and upper bounds for the planar anti-Ramsey numbers when  host graphs are wheels to  obtain   upper bounds for $ar_{_\mathcal{P}}(n, C_6)$  when $n\ge8$ and $ar_{_\mathcal{P}}(n, C_7)$ when $n\geq 13$, respectively.\medskip

We need to introduce more notation.  For a graph $G$ we use $V(G)$, $|G|$,  $E(G)$, $e(G)$, $\delta (G)$ and $\alpha(G)$ to denote the vertex set,  number
of vertices,  edge set, number of edges,   minimum degree, and  independence number of $G$, respectively.
For a vertex $x \in V(G)$, we will use $N_G(x)$ to denote the set of vertices in $G$ which are adjacent to $x$.
We define $N_G[x] = N_G(x) \cup \{x\}$ and $d_G(x) = |N_G(x)|$.
The subgraph of $G$ induced by $A$, denoted $G[A]$, is the graph with vertex set $A$ and edge set $\{xy \in E(G) : x, y \in A\}$. We denote by $B \less A$ the set $B - A$  and $G \less A$ the subgraph of $G$ induced on $V(G) \less A$, respectively.
If $A = \{a\}$, we simply write $B \less a$ and $G \less a$, respectively.
Given two graphs $G$ and $H$,  the \dfn{union} of $G$ and $H$, denoted $G \cup H$, is the graph with vertex set $V(G) \cup V(H)$ and edge set $E(G) \cup E(H)$.
Given two isomorphic graphs $G$ and $H$, we may (with a slight but common abuse of notation) write $G = H$.
Given a plane graph $G$ and an integer $i\ge3$, an \dfn{$i$-face} in $G$ is  a face of size $i$. Let $f_i(G)$ denote the number of $i$-faces in $G$ 
and $n_i(G)$ denote the number of vertices of degree $i$ in $G$. Given an edge-coloring $c$ of $G$, let $c(G)$ denote the number of colors used under  $c$.  For any positive integer $k$, let  $[k]:=\{1,2, \ldots, k\}$. \medskip

\section{Rainbow Paths}\label{path}

In this section, we study  planar anti-Ramsey numbers for paths. We begin with a construction of a plane triangulation $T_H$ that will be needed in the proof of  Theorem~\ref{Pk}.

\begin{lem}\label{wonderful}
For any integers $p\ge1$ and  $n=3p+2$, there exist plane triangulations $H$ on $p+2$ vertices  and $T_H$ on $n$ vertices such that $H$ and $T_H$ satisfy the following. 
\begin{enumerate}[(a)]

\item $H\subseteq T_H$ and $H$ is hamiltonian;
\vspace{-8pt}
\item  $V(T_H) \less V(H)$ 
is an independent set in $T_H$;
\vspace{-8pt}
\item  The   longest  path in $T_H$ has   $2p+5-\max\{0, 3-p\}$ vertices; and
\vspace{-8pt}
\item The longest path in $T_H$  with both endpoints in $V(H)$ has  $2p+3$ vertices.
\end{enumerate}
\end{lem}

\begin{figure}[htbp]
\centering
\includegraphics*[scale=0.3]{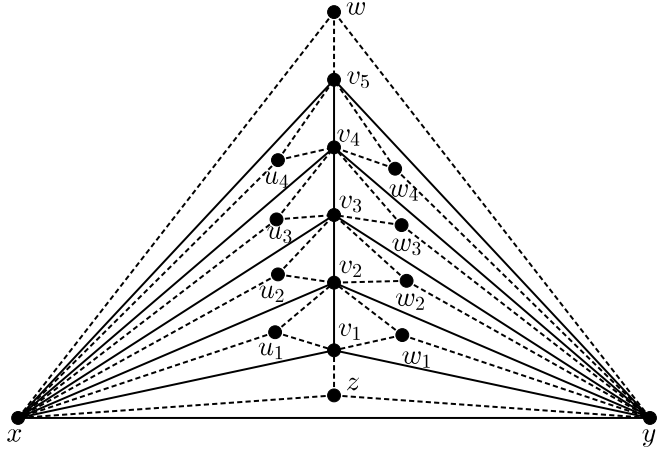}
\caption{The construction of  $T_H$ when $p=5$.}\label{fig1}
\end{figure}

\pf Let $P$ be a path with vertices $v_1, v_2, \dots, v_p$ in order. Let $H$ be the plane triangulation obtained from $P$ by adding two adjacent vertices $x, y$ and joining each of $x$ and $y$ to all vertices on $P$ with the outer face of $H$ having  vertices $x,y,v_p$ on its boundary. Then $|H|=p+2$ and $H$ is hamiltonian. Let $T_H$ be the plane triangulation  obtained from $H$ by adding a new vertex to each $3$-face $F$ of $H$ and   then joining it to all vertices on the boundary of $F$. For each $i\in\{1,2, \dots, p-1\}$, let $u_i$ and $w_i$ be the new vertices added to the faces with vertices $v_i, v_{i+1}, x$ and  $v_i, v_{i+1}, y$ on the boundary, respectively. Let $w, z$ be the new vertices added to the outer face of $H$ and the face of $H$ with vertices $x, y, v_1$ on its boundary. The construction of $T_H$  when $p=5$ is depicted in  Figure \ref{fig1}. Then $|T_H|=|H|+f_3(H)=3p+2=n$ and  $V(T_H) \less V(H)=\{u_1, \dots, u_{p-1}, w_1, \dots, w_{p-1}, w,z\}$. Clearly,  $V(T_H) \less V(H)$ is a maximal independent set of $T_H$ with $|V(T_H) \less V(H)|=f_3(H)=2p$ and  $|V(T_H) \less V(H)|\ge|H|+1-\max\{0, 3-p\}$. 
It can be easily checked that the longest path in $T_H$ has   $2p+5-\max\{0, 3-p\}$ vertices,  and   the  longest path with  both endpoints in $V(H)$ has   $2p+3$ vertices.\qed\medskip

\begin{thm}\label{P8}
For any $ k\in\{8, 9\}$, let  $\varepsilon=k\text{ mod } 2$ and $n\ge k$  be an integer. Then 
$ar_{_\mathcal{P}}(n, P_k)\geq(3n+3\varepsilon-\varepsilon^*-3)/{2}$, where  $\varepsilon^*=(n+1+\varepsilon)\text{ mod }2$.
\end{thm}

\pf   Let $k\in\{8,9\}, n, \varepsilon,\varepsilon^*$ be given as in the statement.
Let $t$ be a positive integer satisfying  $2t-3-\varepsilon+\varepsilon^*=n$. Then $t\geq k-3$ because $n\ge k$.
Let $H$ be a plane drawing of $K_2+(\overline{K}_{t-3-\varepsilon}\,  \cup\,  K_{\varepsilon+1})$. Clearly, $H$ has $3$-faces and $4$-faces only. Notice that  $|H|=t$, $f_3(H)=2+2\varepsilon$, $f_4(H)=t-3-\varepsilon$, $e(H)=2t-3+\varepsilon$, and $H$ is $P_{k-2}$-free but not $P_{k-3}$-free. Let $\mathcal{F}$ be a set which consists of all $4$-faces of $H$ and $\varepsilon^*$ of the $3$-faces of $H$. Let $T^*$ be  the plane triangulation obtained   from  $H$ by adding a new vertex  to each face $F\in \mathcal{F}$ and then joining it to all vertices on $F$.
Then  $|T^*|=|H|+|\mathcal{F}|=|H|+f_4(H)+\varepsilon^*=2t-3-\varepsilon+\varepsilon^*=n$. Clearly, $T^*\in \mathcal{T}_n(P_k)$. Now let $c$ be an edge-coloring of $T^*$ defined as follows:     edges in $E(H)$ are colored with distinct colors under $c$ (that is, $T^*$ contains a rainbow copy of $H$ under $c$), and for each $F\in\mathcal{F}$, all the new edges added inside $F$ are colored the same, but for distinct faces $F, F'\in\mathcal{F}$, new edges inside $F$ are colored differently than the new edges inside $F'$.
It can be easily checked that  $T^*$ has no rainbow $P_k$ but contains a rainbow copy of $P_{k-1}$ under $c$.  Then $c(T^*)=e(H)+f_4(H)+\varepsilon^*=3t-6+\varepsilon^*=(3n+3\varepsilon-\varepsilon^*-3)/2$,  since $n=2t-3-\varepsilon+\varepsilon^*$.
Hence, $ar_{_\mathcal{P}}(n, P_k)\ge c(T^*)\ge(3n+3\varepsilon-\varepsilon^*-3)/2$, as desired. This completes the proof of Theorem~\ref{P8}.\qed\medskip

We next prove a lower bound for  $ar_{_\mathcal{P}}(n, P_k)$ when $k\ge 10$.

\begin{thm}\label{Pk}
Let  $k$ and $n$ be two integers such that  $n\ge k\ge10$. Let $\varepsilon= k\text{ mod }2$.  Then 

\[ar_{_\mathcal{P}}(n, P_k)\geq
\begin{cases}
\,  n+2k-12   & if  \quad  k\le n< 3\lfloor{k}/{2}\rfloor+\varepsilon-5 ,\\[2mm]

\, (3n+9\left\lfloor {k}/{2}\right\rfloor+3\varepsilon-43)/2     & if  \quad 3\lfloor{k}/{2}\rfloor+\varepsilon-5 \le n\le 5\lfloor {k}/{2}\rfloor+\varepsilon-15,\\[2mm]

\, 2n+k-14  & if \quad  n> 5\lfloor {k}/{2}\rfloor+\varepsilon-15.
\end{cases}\]

\end{thm}

\pf  Let $k, n,  \varepsilon$ be given as in the statement.  Assume first  that $ k\le n< 3\lfloor{k}/{2}\rfloor+\varepsilon-5$.  Then $k\ge 12$.  Let $p=k-5$ and let $P$ and  $H$  be  defined in the proof of Lemma \ref{wonderful}. By Lemma~\ref{wonderful},  $|H|=k-3$, $f_3(H)=2k-10$, $e(H)=3k-15$ and $H$ is hamiltonian. Since $n< 3\lfloor{k}/{2}\rfloor+\varepsilon-5$, we see that $n-k+3< f_3(H)$. Let $\mathcal{F}$ be a set which consists of $n-k+3$ many $3$-faces of $H$. Let $T^*$ be the plane triangulation  obtained from $H$ by adding a new vertex  to each face $F\in \mathcal{F}$ and then joining it to all vertices on the boundary of $F$. Clearly, $T^*\in \mathcal{T}_n(P_k)$. Now let $c$ be an edge-coloring of $T^*$ defined as follows: edges in  $E(H)$ are colored with distinct colors  under $c$ (that is, $T^*$ contains a rainbow copy of $H$ under $c$), and for each $F\in\mathcal{F}$, all the new edges added inside $F$ are colored the same, but for distinct faces $F, F'\in\mathcal{F}$, new edges inside $F$ are colored differently than the new edges inside $F'$.
It can be easily checked that $T^*$ has no rainbow $P_k$ but contains a rainbow copy of $P_{k-1}$ under $c$.  Then 
 $c(T^*)=e(H)+|\mathcal{F}|=3k-15+n-k+3=n+2k-12$.
Hence, $ar_{_\mathcal{P}}(n, P_k)\ge c(T^*)\ge n+2k-12$.\medskip

Next assume that $3\lfloor{k}/{2}\rfloor+\varepsilon-5\le  n\le 5\lfloor {k}/{2}\rfloor+\varepsilon-15$.  Let $\varepsilon^*=(n+\lceil k/2 \rceil)\text{ mod }2$.
By the choice of $\varepsilon^*$, let $t$ be a positive integer satisfying $2t+\varepsilon^*+10-3\lfloor k/2 \rfloor-\varepsilon=n$. Since $n\ge 3\lfloor{k}/{2}\rfloor+\varepsilon-5$, it follows that $t-3\lfloor k/2\rfloor+10\ge2+\varepsilon$.
Let $p=\lfloor k/2\rfloor-4$ and let $P$, $H$,  $T_H$, $x, y,  w, v_{\lfloor k/2\rfloor-4}$ be  defined in the proof of Lemma \ref{wonderful}.
By Lemma~\ref{wonderful},   $|H|=\lfloor k/2\rfloor -2$, $f_3(H)=2|H|-4=2\lfloor k/2\rfloor-8$ and $|T_H|=|H|+f_3(H)=3\lfloor k/2\rfloor-10\geq k-5-\varepsilon $.  
Let $F^*$ be the outer face of $T_H$ and $F_0$ be the $3$-face of $T_H$ with vertices $x,w, v_{\lfloor k/2\rfloor-4}$ on its boundary.
Let $T$ be the plane graph on $t$ vertices obtained from $T_H$ by adding $t-3\lfloor k/2\rfloor+10\ge2+\varepsilon$ new vertices to
the face $F^*$ and then joining each of the new vertices to both $x$ and $y$ (and further  adding exactly one    edge among the  new vertices added inside $F^*$ when $ \varepsilon=1$).
Then all $4$-faces of $T$ are inside the face $F^*$ of $T_H$, $e(T)=e(T_H)+2(t-3\lfloor k/2\rfloor+10)+\varepsilon=2t+3\lfloor k/2\rfloor-16+\varepsilon$ and
$f_4(T)=t-3\lfloor k/2\rfloor+10-\varepsilon$.
Let $\mathcal{F}$ be a set which consists of all $4$-faces  of $T$ (and $F_0$ when $\varepsilon^*=1$). Finally, let $T^*$ be the plane triangulation obtained from $T$ by adding a new vertex to each face  $F\in\mathcal{F}$
 and then  joining it to all vertices on the boundary of $F$. Then
$|T^*|=|T|+f_4(T)+\varepsilon^*=2t-3\lfloor k/2\rfloor+10+\varepsilon^*-\varepsilon=n$. By Lemma~\ref{wonderful}, the longest $(x,y)$-path in $T_H$ has  $k-5-\varepsilon$ vertices. Clearly, the longest $(x, y)$-path  in $T^*$ with all its internal vertices inside the face $F^*$   contains all the new vertices  added to $F^*$. Thus  $T^*$ contains   $P_k$ as a  subgraph and so $T^*\in\mathcal{T}_n(P_k)$.   Now  let $c$ be an edge-coloring of $T^*$ defined as follows:  edges in  $E(T)$ are colored with distinct colors under $c$ (that is, $T^*$ contains a rainbow copy of $T$ under $c$),
and for each $F\in \mathcal{F}$, all the new edges added inside $F$ are colored the same, but for distinct $F, F'\in \mathcal{F}$,
new edges inside $F$ are colored differently than the new edges inside $F'$. We see that $T^*$ has no rainbow $P_k$ but contains a rainbow $P_{k-1}$ under $c$.
Since $n=2t+\varepsilon^*+10-3\lfloor k/2\rfloor-\varepsilon$, we see that 
\begin{align*}
c(T^*)=e(T)+f_4(T)+\varepsilon^*
& =(2t+3\left\lfloor{k}/{2}\right\rfloor-16+\varepsilon)
+(t-3\left\lfloor{k}/{2}\right\rfloor+10-\varepsilon)+\varepsilon^*\\
&=(3n+9\left\lfloor {k}/{2}\right\rfloor+3\varepsilon-42-\varepsilon^*)/{2}\\
&\ge (3n+9\left\lfloor {k}/{2}\right\rfloor+3\varepsilon-43)/{2},
\end{align*}
 Hence, $ar_{_\mathcal{P}}(n, P_k)\geq c(T^*)\ge (3n+9\left\lfloor {k}/{2}\right\rfloor+3\varepsilon-43)/2$, as desired.  \medskip

Finally assume that $n\ge5\lfloor {k}/{2}\rfloor+\varepsilon-14$. Let $n-k+7=3m+r$, where $m$ is a positive integer and $r\in\{0,1,2\}$.
Since $k\ge10$ and $n\ge5\lfloor {k}/{2}\rfloor+\varepsilon-14$, we have $m\ge3$ or $m=r=2$.  Let $t:=k+2m-7+\lfloor {r}/2\rfloor$.
Then  $t\ge k-2$ because $m\ge3$ or $m=r=2$, and $t+\lceil(t-k+7)/2\rceil=n-\varepsilon'$, where $\varepsilon'=1$ when $r=1$  and $\varepsilon'=0$ when $r\in\{0,2\}$.  
Let $p=k-9$ and let $P$, $H$,  $x, y,  v_1, \ldots, v_{k-9}$ be  defined in the proof of Lemma \ref{wonderful}.
Then $|H|=k-7$ and  the longest path between $x$ an $y$ in $H$ has $k-7$ vertices.   Let $T'$ be the plane triangulation on $t$ vertices obtained from $H$ by: adding $t-k+7\ge5$ new vertices to the outer face of $H$, then  adding a matching of size $\lfloor (t-k+7)/2\rfloor\ge2$ among  the new vertices,  and finally joining each of the new vertices to both $x$ and $y$.
We see that  $T'$ is a connected $P_{k-2}$-free plane graph with only $3$-faces and $4$-faces. It can be easily checked that  $f_4(T')=\lceil(t-k+7)/2\rceil$ and $e(T')=2t+k-13+\lfloor(t-k+7)/2\rfloor$. Let $F_0$ be the $3$-face  of $T'$ with  vertices $x, y, v_{k-9}$ when $k=10$ and   $x, v_{k-10}, v_{k-9}$ when $k\ge11$ on its boundary. Let $\mathcal{F}$ be a set which consists of  all $4$-faces of $T'$ (and $F_0$ when $\varepsilon'=1$).
 Let $T^*$ be the plane triangulation obtained from $T'$ by adding a new vertex to each $F\in \mathcal{F}$  and then  joining it to all vertices on the boundary of $F$. Then  $|T^*|=|T'|+|\mathcal{F}|=|T'|+f_4(T')+\varepsilon'=t+\lceil(t-k+7)/2\rceil+\varepsilon'=n$. Clearly, $T^*$ contains  $P_k$ as a subgraph and so $T^*\in\mathcal{T}_n(P_k)$. Now let $c$ be an edge-coloring of $T^*$ defined as follows:  edges in $E(T')$ are colored with distinct colors under $c$ (that is, $T^*$ contains a rainbow copy of $T'$ under $c$), and for each
$F\in \mathcal{F}$, all the new edges added inside $F$ are colored the same, but for distinct faces $F, F'\in \mathcal{F}$,
new edges inside $F$ are colored differently than the new edges inside $F'$.  We see that   $T^*$ has no rainbow $P_k$ under $c$ but contains a rainbow copy of $P_{k-1}$. Then 
\begin{align*}
c(T^*)=e(T')+f_4(T')+\varepsilon'& =
\left(2t+k-13+\left\lfloor \frac{t-k+7}{2}\right\rfloor\right)+\left(\left\lceil\frac{t-k+7}{2}\right\rceil+\varepsilon'\right)\\
& =2n+k-13-\varepsilon'+\left\lfloor \frac{t-k+7}{2}\right\rfloor-\left\lceil\frac{t-k+7}{2}\right\rceil\\
&=2n+k-13-\varepsilon'+\left\lfloor \frac{\lfloor {r}/2\rfloor}{2}\right\rfloor-\left\lceil\frac{\lfloor {r}/2\rfloor}{2}\right\rceil\\
&\ge 2n+k-14,
\end{align*}
since $n=t+\lceil(t-k+7)/2\rceil+\varepsilon'$ and $t=k+2m-7+\lfloor {r}/2\rfloor$. Hence, $ar_{_\mathcal{P}}(n, P_k)\geq c(T^*)\ge 2n+k-14$, as desired.  This completes the proof of Theorem~\ref{Pk}.\qed\medskip

\noindent{\bf Remark}. In the proofs of Theorem~\ref{P8} and Theorem~\ref{Pk},   $T^*\in\mathcal{T}_n(P_k)$  has no rainbow $P_k$ but does contain a rainbow copy of $P_{k-1}$ under the coloring $c$ we found. 
\medskip

\section{Rainbow Cycles}\label{cycle}

In this section, we study  planar anti-Ramsey numbers for cycles.

\subsection{Improving  the existing lower bound for $ar_{_\mathcal{P}}(n, C_k)$}
We first prove a lower bound for $ar_{_\mathcal{P}}(n, C_5)$, which improves Theorem~\ref{HJSS2015}(c). \medskip

\begin{thm}\label{C5}
Let  $n\geq 119$ be an integer and let  $ r =(n+7)\text{ mod }18$. Then    $ar_{_\mathcal{P}}(n, C_5)\geq (39n-123-21r)/{18}.$
\end{thm}

\pf   Let $ r, n$ be given as in the statement. Let $t\ge6$ be a positive integer satisfying $18t+11+r=n$. This is possible because   $n\ge119$ and $r= (n+7)\text{ mod }18$.  Let $H$ be  a connected $C_5$-free plane graph with $15t+9$ vertices and  $(12|H|-33)/5$ edges such that  $H$ has only  $3$-faces and $6$-faces, and no  two $6$-faces  share an edge in common. The existence  of such a graph $H$ is due to  Dowden (see  Theorem 4 in  \cite{Dowden}). Notice that $f_6(H)=3t+2$ and $f_3(H)=18t+6$. Let $\mathcal{F}$ be a set which consists of  all $6$-faces  and $r$ of the $3$-faces of $H$. Then $|\mathcal{F}|=f_6(H)+r$.  Let $T^*$ be  the plane triangulation obtained   from  $H$ by adding a new vertex  to each face $F\in \mathcal{F}$ and then joining it to all vertices on the boundary of $F$. Then  $|T^*|=|H|+|\mathcal{F}|=|H|+f_6(H)+r=(15t+9)+(3t+2)+r=18t+11+r=n$ and so $T^*\in\mathcal{T}_n(C_5)$.  Finally let $c$ be an edge-coloring of $T^*$ defined as follows:    edges in $E(H)$ are colored with distinct colors under $c$ (that is, $T^*$ contains a rainbow copy of $H$ under $c$),    and for each $F\in \mathcal{F}$, all the new edges added inside $F$ are colored the same, but for distinct $F, F'\in \mathcal{F}$, new edges inside $F$ are colored differently than the new edges inside $F'$.  We see that $T^*$ has no rainbow $C_5$ under $c$ because $H$ is $C_5$-free and no rainbow $C_5$ in $T^*$ can contain any new edges added to $H$.  Then  $$c(T^*)=e(H)+f_6(H)+r=(36t+15)+(3t+2)+r={(39n-123-21r)}/{18},$$
since $n=18t+11+r$. Therefore, $ar_{_\mathcal{P}}(n, C_5)\ge c(T^*)\ge (39n-123-21r)/{18}$, as desired. This  completes the proof of Theorem~\ref{C5}.\qed\\

\noindent{\bf Remark}.   By Proposition~\ref{LU} and Theorem~\ref{Dowden2015}(c),  $ar_{_\mathcal{P}}(n, C_5)\le ex_{_\mathcal{P}}(n, C_5)\le (12n-33)/{5}$ for  all $n\geq 11$. It then follows from Theorem~\ref{C5}  that $(39n-123-21r)/{18}\le ar_{_\mathcal{P}}(n, C_5) \le ex_{_\mathcal{P}}(n, C_5)\le(12n-33)/{5}$ for all $n\ge119$, where  $ r =(n+7)\text{ mod }18$.\medskip

Theorem~\ref{Ck} below provides a new lower bound for $ar_{_\mathcal{P}}(n, C_k)$ when $k\ge5$, which improves  Theorem~\ref{HJSS2015}(d).

\begin{thm}\label{Ck}
For integers  $k\ge5$, $n\geq k^2-k$, and   $r=(n-2)\text{ mod }(k^2-k-2)$,
$$ar_{_\mathcal{P}}(n, C_k)\geq \left(\frac{k-3}{k-2}+\frac{2}{3(k+1)(k-2)}\right)(3n-6)-\frac{2k^2-5k-5}{k^2-k-2}r.$$
\end{thm}

\begin{figure}[htbp]
\centering
\includegraphics*[scale=0.5]{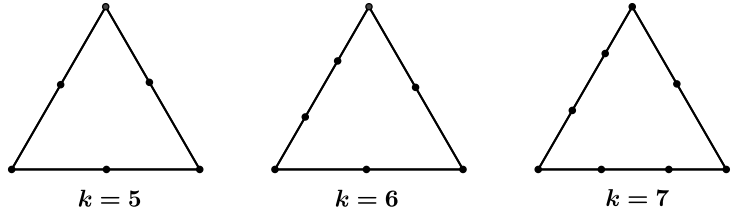}
\caption{Subdividing  one $3$-face  of $T$ when $k\in\{5,6,7\}$.}\label{fig2}
\end{figure}

\begin{figure}[htbp]
\centering
\includegraphics*[scale=0.3]{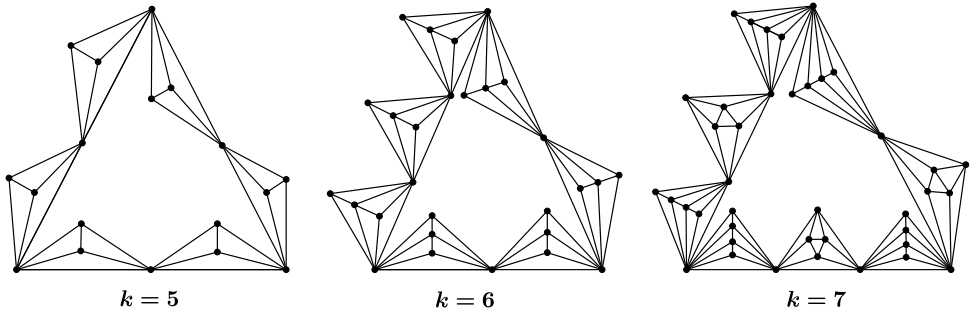}
\caption{ One possible construction of $T'$   when $k\in \{5,6,7\}$.}\label{fig3}
\end{figure}

\pf   Let $n, k, r$ be given as in the statement.   Let   $t\geq 3$ be an integer satisfying $(k^2-k-2)(t-2)+2+r=n$. This is possible because  $r=(n-2)\text{ mod }(k^2-k-2)$  and $n\ge k^2-k$.
Let $T$ be a  plane triangulation on  $t$ vertices. Then  $f_3(T)=2t-4$.  Let $k: =3m+q$, where $q\in\{0,1,2\}$ and $m\ge1$ is an integer.  Let  $T^{\prime}$ be obtained  from $T$ as follows.  For each face $F$ in $T$: 
first subdivide each of the $q+1$ of the edges of $F$  $m$ times;  next,   subdivide each of the remaining $2-q$ edges  of $F$  $m-1$ times;  and finally,   replace each  edge from the subdivision of $T$ by any plane triangulation on $k-1$ vertices.  
Examples of constructions of $T'$ when $k\in\{5,6,7\}$ are depicted  in Figure~\ref{fig2} and Figure~\ref{fig3}.\medskip

 It is worth noting  that different  edges of the subdivision of $T$ may be replaced by different plane triangulations on $k-1$ vertices. Such a subdivision of $T$  is possible when $q\in\{0,1, 2\}$ because when $q=2$, every edge of  $T$  is subdivided $m$ times;  and when  $q\in\{0,1\}$, the dual of $T$ has a perfect matching, say $M$. Let $M^*$ be the dual edges of $M$ in $T$. Then every face $F$ in $T$ contains exactly one edge in $M^*$ and $|M^*|=t-2$.  When $q=0$, each edge in $M^*$ is  divided $m$ times, and when $q=1$, each edge in $M^*$ is divided $m-1$ times.  Thus  $(q+1)(t-2)$ many edges of $T$ are each subdivided $m$ times and $(2-q)(t-2)$ many  edges of $T$ are each subdivided $m-1$ times.  One can check that 
\begin{align*}
|T^{\prime}|&=t+(q+1)(t-2)[(m+1)(k-3)+m]+(2-q)(t-2)[((m-1)+1)(k-3)+ (m-1)]\\
&=t+(t-2)[(q+1)(mk-2m+k-3)+(2-q)(mk-2m-1)]\\
&=t+(t-2)[(q+1+2-q)(mk-2m)+(q+1)(k-3)-(2-q)]\\
&=t+(k^2-k-5)(t-2)
\end{align*}
and 
\begin{align*}
e(T')&=(q+1)(t-2)(m+1)[(3(k-1)-6]+(2-q)(t-2)m[3(k-1)-6]\\
&=(t-2)(3k-9)[(q+1+2-q)m+q+1]\\
&=(t-2)(3k-9)(3m+q+1)\\
&=(t-2)(3k-9)(k+1)=3(k^2-2k-3)(t-2).
\end{align*}

By the construction of $T'$, we see that $T'$ is $C_k$-free (but contains  $C_{k+1}$ as a subgraph),    $T'$ has $f_3(T)$ many $i$-faces with $i>3$ and at least  
\begin{align*}
&(q+1)(t-2)(m+1)[(2(k-1)-5]+(2-q)(t-2)m[(2(k-1)-5]\\
&=(t-2)(2k-7)[(q+1+2-q)m+q+1]\\
&=(t-2)(2k^2-5k-7)\\
&\ge k^2-k-2
\end{align*}
 many $3$-faces because  $t\ge3$ and  $k\ge5$. Let $\mathcal{F}$ be a set which consists of all $i$-faces of $T'$ with $i>4$ and $r$ of the $3$-faces of $T'$.  Let $T^*$ be  the plane triangulation obtained   from  $T'$ by adding a new vertex  to each face $F\in \mathcal{F}$ and then joining it to all vertices on the boundary of $F$.
Then  $|T^*|=|T'|+f_3(T)+r=[t+(k^2-k-5)(t-2)]+(2t-4)+r=(k^2-k-2)(t-2)+2+r=n$ and so $T^*\in \mathcal{T}_n(C_k)$.  Now let $c$ be an edge-coloring of $T^*$ defined as follows:     edges in $E(T')$ are colored  with distinct colors under $c$ (that is, $T^*$ contains a rainbow copy of $T'$ under $c$),    and for each $F\in \mathcal{F}$, all the new edges added inside $F$ are colored the same, but for distinct $F, F'\in \mathcal{F}$, new edges inside $F$ are colored differently than the new edges inside $F'$.  We see that $T^*$ has no rainbow $C_k$  (but contains a rainbow copy of $C_{k+1}$) under $c$ because $T'$ is $C_k$-free  (but contains  $C_{k+1}$ as a subgraph) and no rainbow $C_k$ in $T^*$ can contain any new edges added to $T'$.
Hence,  
\begin{align*}
c(T^*)=e(T^{\prime})+f_3(T)+r
&=3(k^2-2k-3)(t-2)+2(t-2)+r\\
&=(3k^2-6k-7)(t-2)+r\\
&=(3k^2-6k-7)\frac{n-r-2}{k^2-k-2}+r\\
&=\left(\frac{k-3}{k-2}+\frac{2}{3(k+1)(k-2)}\right)(3n-6)-\frac{2k^2-5k-5}{k^2-k-2}r,
\end{align*}
since $n=(k^2-k-2)(t-2)+2+r$. Therefore, $ar_{_\mathcal{P}}(n, C_k)\ge c(T^*)\ge \left(\frac{k-3}{k-2}+\frac{2}{3(k+1)(k-2)}\right)(3n-6)-\frac{2k^2-5k-5}{k^2-k-2}r$. This  completes the proof of Theorem~\ref{Ck}.\qed\medskip

\subsection{New upper bounds for $ar_{_\mathcal{P}}(n, C_k)$ when $k\in\{6,7\}$}

Finally, we use the main  ideas in~\cite{HJSS} to  establish   upper bounds for $ar_{_\mathcal{P}}(n, C_k)$ when $k\in\{6,7\}$. We need to introduce more notation.
Let $C_q$ be a cycle with vertices $v_1, v_2, \dots, v_q$ in order, where $q\ge3$. Let $W_q$ be a wheel obtained from  $C_q$ by adding a new vertex $v$, the \dfn{central vertex}  of $W_q$,  and joining $v$ to  all vertices of  $C_q$.  Vertices $v_1, v_2, \dots, v_q$ are called \dfn{rim vertices} of $W_q$. A cycle $C\subseteq W_q$ is a \dfn{central $k$-cycle} if it contains the central vertex of $W_q$ and $|C|=k$. For any plane triangulation $T$ with at least four vertices and any $v\in V(T)$,  the subgraph of $T$ induced by $N_T[v]$  contains the wheel  $W_{d_T(v)}$ with central vertex $v$ as a subgraph.    Let $c(v)$ be the set of all colors such that each is used to color the edges of $W_{d_T(v)}$  under any edge-coloring $c$ of $T$. 
Lemma~\ref{sumc(v)leqC} below will be used in our proof.

\begin{lem}[\cite{HJSS}]\label{sumc(v)leqC}
Let $T$ be a plane  triangulation and let $c:E(T)\rightarrow [m]$  be a  surjection, where $m$ is a positive integer. Then $$\sum_{v\in V (T)}|c(v)|\geq 4m.$$ \medskip
\end{lem}

To establish an upper bound for $ar_{_\mathcal{P}}(n,C_k)$ when $k\in\{6,7\}$, we use the main ideas in~\cite{HJSS} by studying lower and upper bounds for the planar anti-Ramsey numbers when  host graphs are wheels.   For  integers $k\ge4$ and $q\geq k-1$, we define $ar_{_\mathcal{P}}(W_q,C_k)$ to be  the  maximum number  of colors in an  edge-coloring of $W_q$ that has no rainbow copy of $C_k$. 

\begin{thm}\label{WdCk}
For  integers $k\ge5$ and $q\geq k-1$, $\lfloor\frac{2k-7}{k-3}q\rfloor\leq ar_{_\mathcal{P}}(W_q,C_k) \leq \lfloor\frac{2k-5}{k-2}q\rfloor$.
\end{thm}

\pf  Let $W_q$  be a wheel with rim vertices $v_1, v_2, \dots, v_q$  and central vertex $v$.  To obtain the desired lower bound, let $c:E(W_q)\rightarrow [\lfloor(2k-7)q/(k-3)\rfloor]$ be an edge-coloring of $W_q$ defined as follows: for each  $ i \in [q]$, let $r:=i\text{ mod }(k-3)$ and 
$c(vv_i):=i$,
\[c(v_iv_{i+1})=
\begin{cases}
(k-4)\cdot \frac{i-r}{k-3}+q+r-1, & if \ r\in\{3,4, \dots, k-4\},\\
(k-4)\cdot \frac{i-2}{k-3}+q+1, & if \ r=2,\\
(k-4)\cdot \frac{i}{k-3}+q, & if \ r=0, \\
\end{cases}\]
and
\[c(v_iv_{i+1})=
\begin{cases}
(k-4)\cdot \frac{i-1}{k-3}+q+1, & if \ i\neq q  \ and\ r=1,\\
(k-4)\cdot \frac{i-1}{k-3}+q, & if \ i=q\ and \  r=1, \\
\end{cases}\]
where all arithmetic on the index $i+1$  here and henceforth is done modulo $q$.  It can be easily checked that $c$ is a surjection and $W_q$ has no rainbow $C_k$ (but contains a rainbow copy of $C_{k-1}$) under the coloring $c$.  Hence,  $ar_{_\mathcal{P}}(W_q,C_k)\geq \lfloor(2k-7)q/(k-3)\rfloor$.\medskip

Next we prove that  $ar_{_\mathcal{P}}(W_q,C_k)\leq (2k-5)q/(k-2)$. Let $c:E(W_q)\rightarrow [m]$ be  any  surjection such that $W_q$ contains no rainbow $C_k$ under the coloring $c$. It suffices to show that  $m\le (2k-5)q/(k-2)$.
For any integer $\ell$, let  $A_\ell$ be the set of colors used $\ell$ times under  the coloring $c$. For integers $\alpha\in[m]$ and $j\geq 1$, let:  $\eta_j(\alpha)$ be the number of central $k$-cycles in $W_q$ containing $j$ edges colored $\alpha$ under $c$,  $\eta(\alpha):=\sum_{j=2}^k\eta_j(\alpha)$, $\beta(\alpha):=|
\{i\in [q]: c(vv_{i})=\alpha\}|$ and $\beta^{\prime}(\alpha):=|
\{i\in [q]: c(v_iv_{i+1})=\alpha\}|$. For any integer $\ell$, it is easy to check that $\beta(\alpha)+\beta^{\prime}(\alpha)=\ell$ for any $\alpha\in A_\ell$. Notice that for any integer $i\in[q]$, $vv_i$  belongs to  exactly two central $k$-cycles and $v_iv_{i+1}$ belongs to exactly $k-2$ central $k$-cycles in $W_q$. For any $\alpha\in A_\ell$, we see that $$2\eta(\alpha)\leq 2\eta(\alpha)+\eta_1(\alpha)\leq \sum_{j\geq 1}j\eta_j(\alpha)=2\beta(\alpha)+(k-2)\beta^{\prime}(\alpha)\leq (k-2)\ell,$$
which implies that $\eta(\alpha)\leq (k-2)\ell/2$. Since each of the $q$ central $k$-cycles of $W_q$ contains a color $\alpha$ with $\eta(\alpha)\geq 1$, we have 
$$q\leq \sum_{\ell\geq 2}\sum_{\alpha\in A_\ell}\eta(\alpha)\leq\sum_{\ell\geq 2}(k-2)\ell|A_\ell|/2,$$ which implies $2q/(k-2)\leq \sum_{\ell\geq 2}{\ell|A_\ell|}$. This, together with $2q=e(W_q)=\sum_{\ell\geq 1}{\ell|A_\ell|}$, implies that  $|A_1|\leq (2k-6)q/(k-2)$. Then
\begin{align*}
 m=|A_1|+\sum_{\ell\geq 2}{|A_\ell|}\leq |A_1|+\sum_{\ell\geq 2}|A_\ell|/2
=|A_1|/{2}+\sum_{\ell\geq 1}{\ell}{|A_\ell|}/2
\leq (2k-5)q/(k-2),
\end{align*}
as desired.\qed\medskip

Corollary~\ref{ob} below follows from the fact  that $\lfloor\frac{2k-7}{k-3}q\rfloor=2q-\lfloor\frac{q}{k-3}\rfloor$,  $\lfloor \frac{2k-5}{k-2}q\rfloor=2q-\lfloor\frac{q}{k-2}\rfloor$  and   $ar_{_\mathcal{P}}(W_q,C_k) = 2q-\lfloor \frac{q}{k-3}\rfloor$ if $\lfloor \frac{q}{k-2}\rfloor=\lfloor \frac{q}{k-3}\rfloor$.  One can see that  $\lfloor \frac{q}{k-2}\rfloor=\lfloor \frac{q}{k-3}\rfloor$  when $q\in \{t(k-2), \dots, t(k-2)+k-4-t\}$ for any integer $ t\in [k-4]$.\medskip

\begin{cor}\label{ob}
Let $k\ge5$ and $q\geq k-1$ be integers.  If $q\in \{t(k-2), \dots, t(k-2)+k-4-t\}$ for some integer $ t\in [k-4]$, then $ar_{_\mathcal{P}}(W_q,C_k) = 2q-\lfloor\frac{q}{k-3}\rfloor$.
\end{cor}\medskip

We are ready to determine the exact value for $ar_{_\mathcal{P}}(W_q,C_6)$ when $q\ge5$.\medskip

\begin{thm}\label{WdC6}
For integer $q\geq 5$, $ar_{_\mathcal{P}}(W_q,C_6)=\lfloor{5q}/{3}\rfloor$.
\end{thm}

\pf By Theorem \ref{WdCk},  $ar_{_\mathcal{P}}(W_q,C_6)\geq \lfloor5q/3\rfloor$.
To prove that  $ar_{_\mathcal{P}}(W_q,C_6)\leq \lfloor5q/3\rfloor$, it suffices to show that for any  surjection $c:E(W_q)\rightarrow [m]$ such that  $W_q$  contains no rainbow $C_6$ under the coloring $c$,  we must have $m\le \lfloor5q/3\rfloor$.  We do that next.\medskip

 Let $A_\ell$ be the set of colors used $\ell$ times under the coloring  $c$. For $\alpha\in[m]$, let $\eta_j(\alpha)$ be the number of central $6$-cycles in $W_q$ containing $j$ edges colored $\alpha$ under $c$,  $\eta(\alpha):=\sum_{j=2}^6\eta_j(\alpha)$, $\beta(\alpha):=|
\{i\in [q]: c(vv_{i})=\alpha\}|$ and $\beta^{\prime}(\alpha):=|
\{i\in [q]: c(v_iv_{i+1}) =\alpha\}|$. Then $\beta(\alpha)+\beta^{\prime}(\alpha)=\ell$ for all $\alpha\in A_\ell$. Notice that for any integer $i\in [q]$, $vv_i$ belongs to  exactly two central $6$-cycles and $v_iv_{i+1}$ belongs to exactly  four central $6$-cycles. For any $\alpha\in A_\ell$, we see that $$2\eta(\alpha)\leq 2\eta(\alpha)+\eta_1(\alpha)\leq \sum_{j\geq 1}j\eta_j(\alpha)=2\beta(\alpha)+4\beta^{\prime}(\alpha)\leq 4\ell.$$
This  implies that $\eta(\alpha)\leq 2\ell$. Notice that for  any $\alpha\in A_2$, two edges of $W_q$ colored by $\alpha$ can prevent at most three central $6$-cycles from being rainbow under the coloring $c$, and so $\eta(\alpha)=\eta_2(\alpha)\leq 3$. Since each of the $q$ central $6$-cycles of $W_q$ contains a color, say $\alpha\in[m]$,  with $\eta(\alpha)\geq 1$, it follows that $$q\leq \sum_{\ell\geq 2}\sum_{\alpha\in A_\ell}\eta(\alpha)\leq 3|A_2|+\sum_{\ell\geq 3}{2\ell|A_\ell|}.$$  Thus  $q/2\leq 3|A_2|/2+\sum_{\ell\geq 3}{\ell|A_\ell|}$. This, together with $2q=e(W_q)=\sum_{\ell\geq 1}{\ell|A_\ell|}$, implies that  $2|A_1|+|A_2|\leq 3q$. Then
\begin{align*}
m&=|A_1|+|A_2|+\sum_{\ell\geq 3}{|A_\ell|}\leq |A_1|+|A_2|+\sum_{\ell\geq 3}{\ell}{|A_\ell|}/3\\
&=(2|A_1|+|A_2|)/{3}+\sum_{\ell\geq 1}\ell|A_\ell|/3
=(2|A_1|+|A_2|)/{3}+{2q}/{3}\leq {5q}/{3},
\end{align*}
as desired.\qed\\


Finally,  we obtain new  upper bounds for $ar_{_\mathcal{P}}(n, C_6)$ when $n\geq 8$ and $ar_{_\mathcal{P}}(n, C_7)$ when  $n\geq 13$, respectively.   

\begin{thm}\label{C6}
$ar_{_\mathcal{P}}(n, C_6)\leq {17(n-2)}/6$ for all  $n\geq 8$, and
$ar_{_\mathcal{P}}(n, C_7)\leq {(59n-113)}/{20}$ for all  $n\geq 13$.
\end{thm}

\pf We first prove that   $ar_{_\mathcal{P}}(n,C_6)\leq {17(n-2)}/6$ for  all integers $n\geq 8$.  Let $n\ge8$ be given and  let $T$ be any plane triangulation on $n$ vertices such that $T$ contains $C_6$ as a subgraph. Let $c:E(T)\rightarrow [m]$ be any surjection  such that $T$ contains no rainbow $C_6$ under the coloring $c$.  It suffices to show that  $m\le 17(n-2)/6$. Since $e(T)=3n-6$ and $n\ge8$,  $T$ must have  at least two vertices each with degree at least five. Thus, $n_4(T)\leq n-2-n_3(T)$ and $n_3(T)\ge0$. For any $v\in V(T)$, we see that $|c(v)|\leq e(W_{d_T(v)})=2d_T(v)$. But for any   $v\in V(T)$  with $d_T(v)\geq 5$, by Theorem \ref{WdC6}, $|c(v)|\leq ar_{_\mathcal{P}}(W_{d_T(v)},C_6)=\lfloor5d_T(v)/3\rfloor$. By Lemma \ref{sumc(v)leqC},
\begin{align*}
4m & \leq \sum_{v\in V (T)}|c(v)|\leq 6n_3(T)+8n_4(T)+\sum_{v\in V(T), d_T(v)\geq 5}\left\lfloor5d_T(v)/{3}\right\rfloor\\
&\leq n_3(T)+4n_4(T)/{3}+5/{3}\cdot \sum_{v\in V(T)}d_T(v)\\
& \leq 4(n-2)/{3}-{n_3(T)}/{3}+{5}/{3}\cdot 2(3n-6)\leq 34(n-2)/{3},
\end{align*}
which implies that $m\leq 17(n-2)/6$, as desired.\medskip

It remains to prove that $ar_{_\mathcal{P}}(n,C_7)\leq {(59n-113)}/{20}$ for all $n\geq 13$.  The proof  is similar to the proof of $ar_{_\mathcal{P}}(n,C_6)\leq {17(n-2)}/6$.  We include a proof here for completeness.
 Let $n\ge13$ be given and  let  $T$ be any plane triangulation on $n$ vertices such that $T$ contains $C_7$ as a subgraph. Let $c:E(T)\rightarrow [m]$ be any surjection  such that $T$ contains no rainbow $C_7$ under the coloring $c$. It suffices to show that    $m\le (59n-113)/20$. Since $e(T)=3n-6$ and $n\ge13$,  $T$ must have  at least one vertex of degree six. Thus, $n_5(T)\leq n-1-n_3(T)-n_4(T)$ and $n_i(T)\geq 0$ $(i=3,4)$.  For any $v\in V(T)$, we see that $|c(v)|\leq e(W_{d_T(v)})=2d_T(v)$. But for any   $v\in V(T)$  with $d_T(v)\geq 6$, by Theorem \ref{WdCk},  $|c(v)|\leq ar_{_\mathcal{P}}(W_{d_T(v)},C_7)\leq \lfloor9d_T(v)/5\rfloor$. By Lemma \ref{sumc(v)leqC},
\begin{align*}
4m & \leq \sum_{v\in V (T)}|c(v)|\leq 6n_3(T)+8n_4(T)+10n_5(T)+\sum_{v\in V(T),\ d_T(v)\geq 6}\left\lfloor9d_T(v)/{5}\right\rfloor\\
& \leq {3n_3(T)}/{5}+4n_4(T)/{5}+n_5(T)+{9}/{5}\cdot \sum_{v\in V(T)}d_T(v)\\
&\leq n-1-{2n_3(T)}/{5}-{n_4(T)}/{5}+{9}/{5}\cdot 2(3n-6)\leq {59(n-2)}/{5}+1,
\end{align*}
which implies that $m\leq (59n-113)/20$, as desired. \medskip

This completes the proof of Theorem~\ref{C6}.\qed\\

\noindent{\bf Remark}. A better upper bound for $ar_{_\mathcal{P}}(n, C_6)$  can be obtained using  a result in~\cite{LOSS}  that $ex_{_\mathcal{P}}(n, C_6)\le  18(n-2)/7$ when $n\ge 6$.  By Proposition~\ref{LU} and Theorem~\ref{Ck}, we see that $\frac{65(n-2)}{28}-\frac{37r}{28}\le ar_{_\mathcal{P}}(n, C_6)\le ex_{_\mathcal{P}}(n, C_6)\leq \frac{72(n-2)}{28}$ for all $n\ge30$, where  $r=(n-2)\text{ mod }28$.\\

\noindent {\bf Acknowledgments.} Yongxin Lan and Yongtang Shi are partially supported by National
Natural Science Foundation of China and Natural Science Foundation of Tianjin (No.
17JCQNJC00300).

\end{document}